\documentclass{article}
\usepackage{graphicx,amsmath}
\usepackage{blkarray}
\usepackage{amssymb}
\usepackage{natbib}
\usepackage{todonotes}
\usepackage{mathrsfs}
\usepackage{algorithm} 
\usepackage{multirow}
\usepackage{algpseudocode} 
\usepackage{amsthm}
\usepackage{todonotes}
\usepackage{mathtools} 
\theoremstyle{definition}
\newtheorem{definition}{Definition}[section]
\theoremstyle{theorem}
\newtheorem{theorem}{Theorem}[section]
\newtheorem{corollary}{Corollary}[theorem]  

\newtheorem{proposition}{Proposition}[section]
\theoremstyle{remark}
\newtheorem{remark}{Remark}[section]

\newcommand{\col}[1]{\begin{bmatrix} #1 \end{bmatrix}}
\newtheoremstyle{examplestyle}
  {1em}
  {1em}
  {\normalfont}
  {}
  {\itshape}
  {.}
  {.5em}
  {}

\theoremstyle{examplestyle}
\newtheorem{example}{Example}[section] 
\DeclareMathOperator{\diag}{diag}
\title{State Re-union Maintainability for Semi-Markov Models in Manpower Planning}
\author{Brecht Verbeken and Marie-Anne Guerry}
\date{}

\begin{document}
\maketitle
\begin{abstract}
In previous research the importance of both Markov and semi-Markov models in manpower planning is highlighted. Maintainability of population structures for different types of personnel strategies (i.e. under control by promotion and control by recruitment) were extensively investigated for various types of Markov models (homogeneous as well as non-homogeneous) (\citep{bartholomew1967stochastic}, \citep{vassiliou1984maintainability}).
Semi-Markov models are extensions of Markov models that account for duration of stay in the states. Less attention is paid to the study of maintainability for semi-Markov models. Although, some interesting maintainability results were obtained for non-homogeneous semi-Markov models (\citep{vassiliou1992non}).

The current paper focuses on discrete-time homogeneous semi-Markov models, and explores the concept of maintainable population structures in this setting for a system with constant total size or one with growth factor $1+\alpha$. In particular, a new concept of maintainability is introduced, the so-called State Re-union maintainability ($SR$-maintainability). Moreover, we show that, under certain conditions, the seniority-based paths associated with the $SR$-maintainable structures converge. This allows to characterize the convex set of $SR$-maintainable structures. 
\end{abstract}

\bigskip

\section{Introduction}
Manpower planning is an essential component of contemporary quantitative human resources management. The primary objective of manpower planning is decision support for strategies that address upcoming human resource needs. By doing so, a well-executed manpower planning approach can prevent future deficits or surpluses in staff numbers. This kind of discrepancy between the current and needed workforce is extremely undesirable, as it could result in increased expenses and diminished earnings. In the realm of manpower planning, the use of (semi-)Markov models has gained prominence for their ability to predict future human resource requirements based on current states \citep{de2011markov}. These models offer a probabilistic framework that can capture the transition probabilities between different states of employment over a given time period. After successfully establishing a suitable model, human resources management might aim to attain or maintain a specific workforce composition. This is is where control theory comes into place.

Control theory has a rich and lengthy history in engineering \citep{fernandez2003control}, while control theory for Markov models was first explored in the context of manpower planning in the works of \citep{bartholomew1967stochastic}, \citep{davies1973structural} and \citep{vajda1975mathematical}.
Later on, this research avenue was extended to non-homogeneous Markov chains \citep{vassiliou1984maintainability}, \citep{vasslllou1984stochastic}, \citep{vassiliou1990asymptotically},\citep{georgiou1992periodicity}. Furthermore, for discrete-time homogeneous Markov models the concept of $n$-step maintainability was studied as well \citep{davies1975maintainability}, \citep{guerry1991monotonicity}, \citep{haigh1992further}.

In a manpower Markov model the states correspond to homogeneous personnel groups for which it is assumed that all members have comparable transition probabilities. We call these states, which are often based on the grades within the manpower system, organisational states. The population structure is then defined as the vector that gathers the proportional number of members in each of these organisational states. Control theory deals with the question whether or not a specific desirable population structure can be maintained by an appropriate choice of the parameters that are controllable. 
As to the means of control, the following three categories of flows are considered \citep{bartholomew1967stochastic}:
\begin{itemize}
    \item Wastage, described by the wastage vector $\mathbf{w} = (w_i)$ where $w_i$ is the probability that a member of state $i$ leaves the system within one time interval.
    \item Internal transitions, i.e. promotions and demotions, characterised by the matrix $\mathbf{P}^I$ where $\mathbf{P}^I_{ij}$ is the probability of the internal transition from state $i$ to state $j$ within one time interval.
    \item Recruitment, described by the recruitment vector $\mathbf{r}=(r_i)$ where $r_i$ is the probability that a new members enters the system in state $i$.
\end{itemize}

Although those three categories can be used in order to exert control over the population structure, control by recruitment is often seen as the most attractive means of control. Using wastage as a means of control by dismissing people is frequently seen as unethical and can lead to a decrease in moral and job satisfaction. Whereas adjusting the promotion and demotion rates within an organization might create dissatisfaction among those who feel their expected career progression is being undermined \citep{masood2023uncovering}. Furthermore, it may involve the promotion of inadequately qualified people or the demotion of proficient employees. Hence, control through recruitment is often favoured because it can be implemented without immediate repercussions for those currently employed \citep{bartholomew1967stochastic}.  
However, minor adaptations concerning the means of control, such as the use of the concept of pressure in the states, which was introduced in \citep{kalamatianou1987attainable} or the use of restricted recruitment as in
\citep{ossai2013maintainability},  were considered throughout the literature. Additionally, \citep{georgiou1992partial} introduced the notion of partial maintainability, which entails the need to retain the proportions in the population structure for some states while allowing for flexibility and variability for other states.
In the remainder of this text, we direct our attention to control by recruitment in line with previous work of \citep{davies1973structural} and \citep{vassiliou1984maintainability}. As a consequence, we consider the wastage as well as the internal transition probabilities as fixed and assume that the recruitment vector is under control of management.

Control theory for semi-Markov processes did not get much attention in the literature. Maintainability of the state sizes in the case of a non-homogeneous semi-Markov chains was first studied in \citep{vassiliou1992non}, where Vassiliou and Papadopoulou extended the concept of maintainability by imposing that the number of members is maintained  for each seniority class within an organisational state.
It should be noted that although we follow the historical trend to introduce control theory in Markov processes by using the terminology of manpower planning, this framework might be transferred to other application areas for (semi-) Markov systems.

We will concentrate our attention to the case of discrete-time homogeneous semi-Markov chains, building on the ideas of \citep{davies1973structural} and   \citep{vassiliou1992non}. In Section \ref{def} we introduce a more comprehensive definition of maintainability for semi-Markov models, i.e. the State Re-union maintainability. Afterwards, in Section \ref{region}, we determine the $SR$-maintainable region by uncovering the inherent connection between the State Re-union maintainability of discrete-time homogeneous semi-Markov chains and the (classical) maintainability for non-homogeneous Markov chains. Finally, in Section \ref{section5}, we make a comparison between the $SR$-maintainable region and the (classical) maintainable region for a Markov chain.

The following section provides the essential background information concerning maintainability under control by recruitment for a Markov chain.

\section{Maintainability under control by recruitment for a Markov chain}

Let us consider a Markov system with $l$ organisational states $S_1, ..., S_l$. The stock vector of the system at time $t$ is the vector with $i$-th element the number of members in  $S_i$ at time $t$. By expressing the stocks proportionally, one obtains a  probability vector $\mathbf{s}(t)$ that is called the population structure of the system at time $t$. Denoting the internal transition matrix as $\mathbf{P}_M^I$, where $M$ stands for Markov, the wastage vector as $\mathbf{w}$ and the recruitment vector, which is a probability vector, as $\mathbf{r}$, the evolution of the population structure for a system with constant total size, can be described as in \citep{bartholomew1967stochastic}:
\begin{equation} \label{Pwr}
    \mathbf{s}(t+1) = \mathbf{s}(t) (\mathbf{P}_M^I + \mathbf{w'} \mathbf{r} )
\end{equation}
and a structure $\mathbf{s}$ is maintainable under control by recruitment if and only if
\begin{equation*} \label{maint}
    \mathbf{s} \geq \mathbf{s} \, \mathbf{P}_M^I
\end{equation*}
It is important to emphasize that the internal transition matrix $\mathbf{P}_M^I$ as we consider here, is a substochastic matrix. It becomes a probability matrix by extending the set of organisational states, which are homogeneous personnel groups within the system, with an additional absorbing state: the wastage state. This extended stochastic matrix is denoted by $\mathbf{P}_M^S$.
\begin{definition} \label{rplus}
    For a maintainable structure $\mathbf{s}$ we define the non-normalised recruitment vector $\mathbf{r^{+}}$ as follows:
    \[\mathbf{r^{+}}:= \mathbf{s} - \mathbf{s} \, \mathbf{P}_M^I\]
\end{definition}

By normalising the vector  $\mathbf{r^+}$ with respect to the $L^1$ norm, one gets the recruitment vector $\mathbf{r}$ that maintains the population structure $\mathbf{s}$.

Furthermore, in \citep{bartholomew1967stochastic} the set of maintainable structures is fully described as the convex hull of its vertices:
\begin{theorem} \citep{bartholomew1967stochastic} \label{the:mainmarkov}
Suppose that $\mathbf{P}_M^I$ is the internal transition matrix of a Markov system. The maintainable region for a system with constant total size and under control by recruitment is the convex set with vertices given by the normalised rows of $(\mathbf{I}- \mathbf{P}_M^I)^{-1}$, with respect to the $L^1$ norm.
\end{theorem}
\begin{remark}
    The reader should be aware that in the remainder of the text, whenever we speak of normalised rows, we are referring to the process of forming probability vectors using the $L^1$
  norm on specific vectors or matrix rows.
\end{remark}
The matrix $\mathbf{I}$ is the $l \times l$ identity matrix. Since $ \mathbf{P}_M^I$ is a substochastic matrix the series $\mathbf{I} +  \mathbf{P}_M^I +  \left(\mathbf{P}_M^I\right){^2} + ...$ converges to $(\mathbf{I}- \mathbf{P}_M^I)^{-1}$. Which not only ensures that all entries of $(\mathbf{I} - \mathbf{P}_M^I)^{-1}$ are non-negative but guarantees the invertibility of $\mathbf{I}- \mathbf{P}_M^I$. 
Hence, each row of $(\mathbf{I}- \mathbf{P}_M^I)^{-1}$ is a non-negative vector that, after normalisation, results in a probability vector. Those vectors span the set of maintainable structures. We introduce the following notion:
\begin{definition}( $(k-1)$ Probability simplex) \\
The $(k-1)$-dimensional probability simplex 
$\Delta^{k-1}$ 
is the simplex whose vertices
are the $k$ standard unit vectors in 
$\mathbb{R}^{k}$, or in other words:
\begin{equation*}
     \Delta^{k-1} = \left\{ (x_1, x_2, \ldots, x_k) \in \mathbb{R}^k \,\middle|\, x_i \geq 0, \sum_{i=1}^k x_i = 1 \right\}
\end{equation*}
\end{definition}

We further generalise Theorem \ref{the:mainmarkov} to systems with growth factor $1+\alpha$ \, (for $\alpha \geq -1$), i.e. where $\mathbf{N}(t+1)$, the total number of people in the system at time $t+1$ equals $(1+\alpha)\mathbf{N}(t)$. In order to do so  the notion of $M$-matrices is useful and we therefore revise some important and known properties of $M$-matrices.\newpage

\begin{definition}($M$-matrix) \citep{berman1994nonnegative}\\
A matrix of the form $\mathbf{A}=d\mathbf{I}-\mathbf{B}$, where $\mathbf{B}$ is a non-negative matrix and $d \geq \rho(\mathbf{B})$, is called an $M$-matrix.
\end{definition} 

Hereby the notations $\rho(\mathbf{B}) = \max_{\lambda \in \sigma(\mathbf{B})} |\lambda|$ and $\sigma(\mathbf{B})$ respectively refer to the spectral radius and the spectrum of the matrix $\mathbf{B}$. Furthermore, remark that if $d=\rho(\mathbf{B})$, it follows that $d$ is the Perron-Frobenius eigenvalue of $\mathbf{B}$, which implies the singularity of $d\mathbf{I}-\mathbf{B}$. We conclude: 

\begin{proposition} \label{nsm}
    An $M$-matrix $\mathbf{A}=d\mathbf{I}-\mathbf{B}$ is invertible if and only if  $d > \rho(\mathbf{B})$.
\end{proposition}

\begin{proposition} \label{non-singular M-matrix} \citep[Theorem 6.2.3]{berman1994nonnegative} For a non-singular matrix $\mathbf{A}$, the following are equivalent:
\begin{itemize}
    \item $\mathbf{A}$ is an $M$-matrix
    \item $\mathbf{A}$ is inverse-positive, i.e. $\mathbf{A}^{-1}$ exists and is a non-negative matrix.
\end{itemize}
    
\end{proposition}
The procedure to determine the maintainable region of a Markov system with internal transition matrix $\mathbf{P}_M^I$, as introduced in \citep{bartholomew1967stochastic} and resulting in Theorem \ref{the:mainmarkov}, requires that the inverse of  $\mathbf{I}- \mathbf{P}_M^I$ exists and is non-negative. We will extend this reasoning to the case of a system under control by recruitment with growth factor $1+\alpha$. Remark that a negative $\alpha$ leads to a contracting system, while a positive $\alpha$ corresponds to a expanding system and a constant sized system coincides with $\alpha=0$.

This leads to the notion of the generalised maintainability equation \citep{bartholomew1967stochastic}:
\begin{definition}(Maintainability equation) \\
A structure $\mathbf{s}$ is maintainable for a system with growth factor $1+\alpha$ under control by recruitment, if the following equation holds for a
non-negative vector $\mathbf{r}^+$:
\begin{equation} \label{maineq2}
    (1+\alpha) \mathbf{s} = \mathbf{s}\mathbf{P}_M^I+\mathbf{r}^+
\end{equation}
If $\mathbf{r}^+ \not  = 0$, we can normalise this vector to get the probability vector $\mathbf{r}$. We then can reformulate the  maintainability equation using  $\mathbf{r}$ as well:
\begin{equation} \label{maineq}
    (1+\alpha) \mathbf{s} = \mathbf{s}\mathbf{P}_M^I+(\mathbf{s} \mathbf{w}^{'}+ \alpha) \mathbf{r}
\end{equation}
    
\end{definition}
It is often more convenient to study the stock vectors instead of the structures $\mathbf{s}$. Writing $\mathbf{n}(t)$ for the stock vector at time $t$, we can say that a structure $\mathbf{s}= \frac{\mathbf{n}(0)}{||\mathbf{n}(0) ||_1}$ is maintainable if there exists a path of stock vectors $(\mathbf{n}(t))_t$ for which a sequence of non-normalised recruitment vectors $(\mathbf{r}^+(t))_t$ exists, such that, for all $t \in \mathbb{N}$:
\begin{equation} \label{maintainable stock vector}
    \mathbf{n}(t+1)=\mathbf{n}(t) \mathbf{P}_M^I+\mathbf{r}^+(t)=(1+\alpha)\mathbf{n}(t)
\end{equation}
Indeed, Equation (\ref{maintainable stock vector}) expresses a necessary and sufficient condition to have $\frac{\mathbf{n}(t)}{||\mathbf{n}(t) ||_1}=\mathbf{s}$,  for all $t \in \mathbb{N}$. 
\newpage
Remark that, due to the strict substochasticity of $\mathbf{P}_M^I$, the number of recruitments is strictly positive in the case of a constant size system or a system subjected to growth ($\alpha \geq 0$). However zero-recruitment might be a possibility for a contracting system ($-1 \leq \alpha < 0$). \\
In case of zero-recruitment we know that $\mathbf{r}^+=0$. Hence, we get that
\begin{equation*}
    \mathbf{s} \mathbf{P}_M^I +\mathbf{r}^+ =  \mathbf{s} \mathbf{P}_M^I=(1+\alpha) \mathbf{s}
\end{equation*}
So, a maintainable population structure $\mathbf{s}$ is necessarily an eigenvector of  $\mathbf{P}_M^I$ with eigenvalue $(1+\alpha)$. 
\begin{example}
    Suppose that $\mathbf{P}_M^I=\begin{pmatrix}
        0.5 & 0.4 & 0 \\
        0 & 0.6 & 0.3 \\
        0 & 0 & 0.8
    \end{pmatrix}$ and $\alpha=-0.2$, so \\ $1+\alpha=0.8$. An easy calculation shows that $\mathbf{s}=(0,0,1)$ is an eigenvector of  $\mathbf{P}_M^I$ with eigenvalue $0.8$, and, as $0.8$ is the growth factor of this system, it follows that the structure $\mathbf{s}$ is a  maintainable structure. 
\end{example}
Moreover, remark that the possibilities for a contracting system are limited  by the wastage vector $\mathbf{w}$. The elements of 
$\mathbf{w}$ set an upper limit on the fraction of people that can exit from each state. For any given state, the system can not possibly lose more individuals than the proportion specified by its corresponding element in $\mathbf{w}$. I.e.
$\alpha < - \max {w}_i$ is not possible. In addition, the following proposition could come to use when one wants to further narrow the possible values for $\alpha$:
\begin{proposition} \label{a en nonzero column}
    Suppose that the $i$th column of $\mathbf{P}_M^I$ possesses no zeros. If the structure $\mathbf{s}$ is maintainable for a system with growth factor $1+\alpha$, then 
    \begin{equation} \label{restriction on alpha}
        \alpha \geq (\mathbf{P}_M^I)_{ii} - 1
    \end{equation}
\end{proposition}
\begin{proof}
    For a maintainable structure $\mathbf{s}$, one obtains that 
    \begin{equation*}
        s_i ({\mathbf{P}_M^I})_{ii} \leq (1+\alpha) s_i
    \end{equation*}
    We claim that $s_i=0$ is not possible. By contradiction, suppose that $s_i=0$, as $\mathbf{s}$ is a probability vector, it possesses a non-zero element, say $s_j$. Furthermore $s_j  ({\mathbf{P}_M^I})_{ji} $ is non-zero, as $({\mathbf{P}_M^I})_{ji}$ is non-zero. But then it follows that $(\mathbf{s}{\mathbf{P}_M^I})_{i}$ is non-zero, while $s_i=0$. Resulting in a contradiction with $\mathbf{s}$ being maintainable. We conclude that
    \begin{equation*}
        ({\mathbf{P}_M^I})_{ii} \leq 1+\alpha 
    \end{equation*}
    and Equation (\ref{restriction on alpha}) follows.
\end{proof}
This proposition can be rephrased into the following corollary:
\begin{corollary}
    Suppose that the $i$th column of $\mathbf{P}_M^I$ possesses no zeros. If  
    \begin{equation*}
        ({\mathbf{P}_M^I})_{ii} > 1+\alpha 
    \end{equation*}
    then the maintainable region for the corresponding Markov system with growth factor $1+\alpha$ is empty.
\end{corollary}
Values of $\alpha$ that satisfy $\alpha < - \max {w}_i$ or that are not in accordance with Proposition \ref{a en nonzero column}, result in an empty maintainable region for the corresponding Markov system.
A complete characterisation of the maintainable region is given in the following theorem:

\begin{theorem} (Characterisation of the Maintainable Region) \label{markovexco}
For a system with internal transition matrix $\mathbf{P}_M^I$ and growth factor $1+\alpha$ the maintainable region satisfies the following:
    \begin{enumerate}
    \item  The maintainable region is a convex set that can be determined as the intersection of the $l$ halfspaces given by:
    \[ (\mathbf{s} \mathbf{P}_M^I)_i \leq (1+\alpha) \mathbf{s}_i \]

    \item If \(  \mathbf{(1+\alpha)I}-\mathbf{P}_M^I \) is invertible then the maintainable region consists of the  convex combinations of the normalised rows of $(\mathbf{(1+\alpha)I}-\mathbf{P}_M^I)^{-1}$, that belong to $\Delta^{l-1}$.

    \item If \( \alpha > \rho(\mathbf{P}_M^I) - 1 \) then the vertices of the maintainable region correspond to the normalised rows of \( (\mathbf{(1+\alpha)I}-\mathbf{P}_M^I)^{-1} \).
\end{enumerate}
\end{theorem}

\begin{proof}
Being an element of the intersection of the halfspaces $(\mathbf{s} \mathbf{P}_M^I)_i \leq (1+\alpha) \mathbf{s}_i$ is equivalent with the existence of a non-normalised recruitment vector $\mathbf{r}^+$. Additionally, rearranging the maintainability Equation (\ref{maineq}) results in:
    \begin{equation*}
        \mathbf{s} \Big((1+\alpha) \mathbf{I}-\mathbf{P}_M^I\Big)= \mathbf{s} \mathbf{w}^{'} \mathbf{r} + \alpha \mathbf{r}
    \end{equation*}
    If $(1+\alpha) \mathbf{I}-\mathbf{P}_M^I$ is invertible we can follow the procedure that resulted in Theorem \ref{the:mainmarkov} by replacing $\mathbf{I}-\mathbf{P}_M^I$ by $(1+\alpha) \mathbf{I}-\mathbf{P}_M^I$. Hence, the maintainable structures $\mathbf{s}$ can be expressed as convex combinations of the normalised rows of $(\mathbf{(1+\alpha)I}-\mathbf{P}_M^I)^{-1}$. Observe that only vectors in the probability simplex $\Delta^{l-1}$  are allowed, so only the convex combinations of the normalised rows of $(\mathbf{(1+\alpha)I}-\mathbf{P}_M^I)^{-1}$ that belong to $\Delta^{l-1}$, should be taken into account as elements of the maintainable region.\\ If  $\alpha > \rho(\mathbf{P}_M^I) - 1$  then $\mathbf{(1+\alpha)I}-\mathbf{P}_M^I$ is a non-singular $M$-matrix by Proposition \ref{nsm}. Hence, according to Proposition \ref{non-singular M-matrix}, $(1+\alpha) \mathbf{I}-\mathbf{P}_M^I$ is inverse-positive, which implies that all of the convex combinations of the normalised rows of $\mathbf{(1+\alpha)I}-\mathbf{P}_M^I$ belong to $\Delta^{l-1}$.
\end{proof}

\section{SR-maintainability under control by recruitment for a semi-Markov chain} \label{def} 
When using semi-Markov chains, the main goal is often to estimate the semi-Markov kernel $\mathbf{q}$ \citep{barbu2009semi}. Suppose that the organisational states of a semi-Markov chain are given by $\mathbf{S}=\{S_1, S_2, \cdots, S_l\}$ and let $T_n$ and $J_n$ denote the time of the $n$-th transition and the state occupied after the $n$-th transition, respectively. The semi-Markov kernel ${\mathbf{q}} = (q_{ij}(k) :  1 \leq i, j  \leq l, k \in \mathbb{N})$ where 
\[ q_{ij}(k) = \Pr(J_{n+1} = S_j, T_{n+1} - T_n = k | J_n = S_i). \]

This semi-Markov kernel can be used to construct a sequence of matrices $\{\mathbf{P}(k)\}_{0 \leq k \leq K}$ where $\mathbf{P}(k)$ corresponds to the one-step ahead transition matrix for staff with organisational state seniority $k$. In practice, $k$ is limited to $\{0, ..., K\}$, where $K$ denotes the maximal attainable state seniority, i.e. the maximal duration of stay in an organisational state.

\begin{theorem} \label{bigmatrix} \citep{verbeken2021discrete} For all $k$ such that \\ $\sum_{h\in\mathbf{S}}\sum_{m=0}^{k-1}q_{ih}\left(m\right)\neq1$ we have:
$$\mathbf{P}_{ij}(k)=
\begin{cases}
    \dfrac{q_{ij}\left(k\right)}{1-\sum_{h\in\mathbf{S}}\sum_{m=0}^{k-1}q_{ih}\left(m\right)} & \text{if } i \not = j \\
   1- \sum_{i \not = j} \dfrac{q_{ij}\left(k\right)}{1-\sum_{h\in\mathbf{S}}\sum_{m=0}^{k-1}q_{ih}\left(m\right)}  & \text{if }  i = j
\end{cases}$$
\end{theorem}
\begin{remark}
    Following the classical notation \citep{barbu2009semi}, where longitudinal data over the time interval $[0,M]$ is used, writing
    $N_{ij}(k)$ for the number of persons in state $S_i$ with organisational state seniority $k$ that go to state $S_j$ during the next time step, while writing $N_i(k)$ for the total number of people in state $S_i$ with organisational state seniority $k$ over the time interval $[0, M-1]$ ,  
    the following maximum likelihood estimator can be used for $\mathbf{P}_{ij}(k)$:
    \begin{equation}\label{schatting P_{ij}(k)}
        \widehat{\mathbf{P}_{ij}(k)}=\frac{N_{ij}(k)}{N_i(k)}
    \end{equation}
\end{remark}
First of all, in order to study the maintainability of a semi-Markov chain, we need to combine all of the information in the matrices $\mathbf{P}(k)$ into one matrix $\mathbf{P}_{SM}$ that characterises the semi-Markov (SM) model.
This can be done by breaking up the states $\mathbf{S}$ according to organisational state seniority. The matrix $\mathbf{P}_{SM}$ is then the transition matrix in accordance with these seniority based disaggregated states. 
\begin{definition} \label{general_P} Suppose we have $l$ organisational states and one state that corresponds to leaving the system, which is called the wastage state. If the sequence $\{\mathbf{P}(k)\}_k$ is of length $K+1$, then for $0\leq k\leq K$ the elements of $\mathbf{P}_{SM}$ equal:
\begin{equation*}
    \left(\mathbf{P}_{{SM}} \right)_{ij}=0 \quad {\text{for } i-1  \not \equiv_{K+1}  k  }
\end{equation*}
and, if $ i-1 \equiv_{K+1}  k $:
\begin{equation} \label{bigm}
\left(\mathbf{P}_{{SM}} \right)_{ij}=
\left\{
    \begin{array}{lll}
    \mathbf{P}(k)_{\lceil \frac{i}{K+1} \rceil, \lceil \frac{i}{K+1} \rceil} & 
    \text{if } \lceil \frac{i}{K+1} \rceil = \lceil \frac{j}{K+1} \rceil & 
    \text{ and } (j-1-i) \equiv_{K+1} 0 \\
    \mathbf{P}(k)_{\lceil \frac{i}{K+1} \rceil, \lceil \frac{j}{K+1} \rceil} & 
    \text{if } \lceil \frac{i}{K+1} \rceil \neq \lceil \frac{j}{K+1} \rceil & 
    \text{ and } (j-1) \equiv_{K+1} 0 
\end{array}
\right.
\end{equation}
\end{definition}
Where we used the notation $a \equiv_{n}  b$ for $a=b \mod n$.
In general, if we describe a transition from $S_i$ to $S_j$, the number $(i-1) \mod (K+1)$  corresponds to the organisational state seniority $k$ in the organisational state where the transition starts from, while $\lceil \frac{i}{K+1} \rceil$ and $\lceil \frac{j}{K+1} \rceil$ correspond to the organisational states themselves. 
This implies that Equation (\ref{bigm}) expresses the two possible types of internal transitions: Either, the organisational state is preserved, which consequently means that the organisational state seniority increases by one, or the organisational state changes, which implies that the organisational state seniority is reset to zero.  
\begin{example}
Suppose that we have a sequence of $4 \times 4$ one-step ahead transition matrices $\{\mathbf{P}(k)\}_{0 \leq k \leq 3}$ of length $4$ with $2$ internal states and one wastage state, so $l=2$. Then according to (\ref{bigm}) it follows that 
\begin{equation*}
\mathbf{P}_{SM}=
\begin{bmatrix}
0 & P_{11}(0) & 0 & 0 &
P_{12}(0) & 0 & 0 & 0 \\
0 & 0 & P_{11}(1) & 0 &
P_{12}(1) & 0 & 0 & 0 \\
0 & 0 & 0 & P_{11}(2) &
P_{12}(2) & 0 & 0 & 0 \\
0 & 0 & 0 & 0 &
P_{12}(3) & 0 & 0 & 0 \\
P_{21}(0) & 0 & 0 & 0 &
0 & P_{22}(0) & 0 & 0 \\
P_{21}(1) & 0 & 0 & 0 &
0 & 0 & P_{22}(1) & 0 \\
P_{21}(2) & 0 & 0 & 0 &
0 & 0 & 0 & P_{22}(2) \\
P_{21}(3) & 0 & 0 & 0 &
0 & 0 & 0 & 0 \\
\end{bmatrix}
\end{equation*}
\end{example}
\begin{remark}
    In line with the existing literature on maintainability \citep{bartholomew1967stochastic}, we will assume that for each seniority-based state the wastage probability is non-zero, which implies that the matrix $\mathbf{P}_{SM}$ is strictly substochastic, i.e. all row sums of $\mathbf{P}_{SM}$ are strictly smaller than $1$.
\end{remark}
Redefining the state set may be more practical in this context. Therefore, we define the following:
\begin{definition}
The set of seniority-based states is given by $$\mathbf{S_{SB}}=\{S_{a(b)} \mid 0 \leq a \leq K \text{ and } 1 \leq b \leq l \},$$ where the state $S_{a(b)}$ corresponds to the staff in organisational state $b$ that has  organisational state seniority equal to $a$.
\end{definition}
This allows us to consider $\mathbf{P}_{SM}$ as the transition matrix with state space $\mathbf{S_{SB}}$. As is the case for the Markov system, $\mathbf{P}_{SM}$ is no probability matrix but results into one after supplementing the set of seniority-based states with the wastage state and its corresponding probabilities.

 By writing for the state space $\mathbf{S_{SB}}$, the stock vector at time $t$ by ${\mathbf{n}}_{\mathbf{SB}}(t)$ and the non-normalised recruitment vector at time $t$ by ${\mathbf{r}}_{\mathbf{SB}}^{+}(t)$, we obtain the following equations, that describe the evolution of the stock vector for a system with growth factor $1+\alpha$: \newpage
\begin{equation}\label{SR1bis}
   {\mathbf{n}}_{\mathbf{SB}}(t+1)= {\mathbf{n}}_{\mathbf{SB}}(t) \mathbf{P}_{SM}+    {\mathbf{r}}_{\mathbf{SB}}^{+}(t)
\end{equation}
\begin{equation}
    \mathbf{N}(t+1)=(1+\alpha)\mathbf{N}(t)
\end{equation}

Now the question is, what should be maintained? Which proportions should be preserved? In the case of a Markov chain, often, a fixed proportion for each of the states is to be preserved. We could do the same in the semi-Markov case, but this would mean that the proportions of all seniority-based states should be preserved. This is in line with the approach in \cite{vassiliou1992non}, where this notion of maintainability was introduced for non-homogeneous semi-Markov chains. However, in the context of real-world applications, it might often be deemed less restrictive and more  efficacious to focus on preserving the proportions in the organisational states instead.
In other words, for each 
$b \in \mathbf{S}$ 
we fuse the states $\{S_{a(b)}\}_{0 \leq a \leq K}$ into one state.
This leads to a generalised definition of maintainability for semi-Markov chains:

\begin{definition}(Re-union matrix) \\
For a transition matrix $\mathbf{P}_{SM}$ with state space $\mathbf{S_{SB}}$, a $(K+1)l \times l$ matrix $\mathbf{U}=(U_{ij})$ is called the re-union matrix if each of its $l$ columns $\col{\mathbf{U}}_j$ consists of $K+1$ ones through:
\begin{equation*}
    U_{ij}=
    \begin{cases}
    1 & \text{if } (j-1)(K+1) \leq i \leq j(K+1) \\
    0 & \text{ else }
\end{cases} 
\end{equation*}
\end{definition}

If we combine this idea of a re-union matrix with Equation (\ref{SR1bis}), we arrive at the main concept of this paper, being State Re-union maintainability.

\begin{definition}(State Re-union maintainability) \label{M-state}\\
A structure ${\mathbf{s}} \in  \Delta^{l-1}$ is called State Re-union maintainable ($SR$-maintainable) for a system with growth factor $1+\alpha$ under control by recruitment if there exists a path of seniority-based stock vectors 
    $(\mathbf{n}_{\mathbf{SB}}(t))_t$ 
    and if a sequence of recruitment vectors $(\mathbf{r}_{\mathbf{SB}}^{+}(t))_t$  can be chosen
      such that, for every $t \in \mathbb{N}$:
    \begin{equation} \label{SR1}
     {\mathbf{n}}_{\mathbf{SB}}(t+1)= {\mathbf{n}}_{\mathbf{SB}}(t) \mathbf{P}_{SM}+    {\mathbf{r}}_{\mathbf{SB}}^{+}(t)
\end{equation}
\begin{equation} \label{SR2}
    (1+\alpha){\mathbf{n}}_{\mathbf{SB}}(t)\mathbf{U}={\mathbf{n}}_{\mathbf{SB}}(t+1)\mathbf{U}
\end{equation}
\begin{equation} \label{SR3}
\mathbf{s}=\frac{\mathbf{n}_{\mathbf{SB}}(t) \mathbf{U}}{||\mathbf{n}_{\mathbf{SB}}(t) \mathbf{U}||_1} 
\end{equation}
A sequence of seniority-based stock vectors 
$(\mathbf{n}_{\mathbf{SB}}(t))_t$ that satisfies equations (\ref{SR1},\ref{SR2}) and (\ref{SR3}) will be called a seniority-based path associated to the $SR$-maintainable population structure ${\mathbf{s}}$.

\end{definition}

\begin{remark}
    Recruitment is only allowed in seniority-based states with zero seniority. So the recruitment vectors ${\mathbf{r}}_{\mathbf{SB}}^{+}(t)$ can solely possess non-zero  entries for the seniority-based states with zero seniority.
\end{remark}

\section{The $SR$-maintainable region} \label{region}
Now that we have defined the concept of State Re-union maintainability, we will focus our attention to the $SR$-maintainable region $\mathscr{MR}_{SM}(1+\alpha)$ which is the set of all State Re-union maintainable structures for a semi-Markov system with growth factor $1+\alpha$.  
\begin{definition} ($SR$-maintainable region) \\
The $SR$-maintainable region $\mathscr{MR}_{SM}(1+\alpha)$ is the set of State Re-union maintainable structures $\mathbf{s}$ for a semi-Markov system with growth factor $1+\alpha$.
\end{definition}
We characterise and describe this region for a constant sized system in Section \ref{SRMcon} and further generalise this for $\alpha \neq 0$ in Section \ref{SRMo}. A concrete example is provided in Section \ref{section5} where the $SR$-maintainable region and maintainable regions for induced Markov chains are compared.
\subsection{The $SR$-maintainable region for a system with constant size} \label{SRMcon}
As is customary in the study of maintainable regions in the classical sense, we will study the normalised version of the stock vectors, i.e. the population structures \begin{equation*}
   {\mathbf{s}}_{\mathbf{SB}}(t):=\frac{{\mathbf{n}}_{\mathbf{SB}}(t)}{||{\mathbf{n}}_{\mathbf{SB}}(t)||_1}.
\end{equation*}
In general one might note that we already know a subset of the $SR$-maintainable region, namely the subset of population structures for which a seniority-based path does exist with ${\mathbf{s}}_{\mathbf{SB}}(t)={\mathbf{s}}_{\mathbf{SB}}$ i.e. where ${\mathbf{s}}_{\mathbf{SB}}$ is a maintainable vector in the traditional sense for the matrix $\mathbf{P}_{SM}$. What might be more surprising is the following theorem, which shows that all seniority-based paths associated to a $SR$-maintainable population structure converge to such a ${\mathbf{s}}_{\mathbf{SB}}$ that is maintainable for $\mathbf{P}_{SM}$. \\

\begin{theorem} \label{constant}
    Suppose that we have a matrix $\mathbf{P}_{SM}$ which characterises a semi-Markov model, and the re-union matrix $ \mathbf{U}$. Denote the maximal organisational state seniority by $K$. If the structure $\mathbf{s}$ is $SR$- maintainable, then, it follows, for a seniority-based path $\mathbf{s}_{\mathbf{SB}}(t)$ associated to $\mathbf{s}$ that the sequence $\mathbf{s}_{\mathbf{SB}}(t)$, converges, i.e
    \begin{equation}\label{state based structucture converges}
        \lim_{t \to \infty} 
        \mathbf{s}_{\mathbf{SB}}(t) = \mathbf{s}_{\mathbf{SB}}^*
    \end{equation}
\end{theorem}
\begin{proof}
We will prove the convergence of $\mathbf{s}_{\mathbf{SB}}(t)$ by demonstrating the convergence of its constituent parts. More particularly, we partition the set $S_{\mathbf{SB}}$ of seniority-based states according to the organisational state where they are related to. For the organisational state $S_i$ we  denote the set of related seniority-based states as $S_{\mathbf{SB}}[i]= \{S_{a(i)}| 0 \leq a \leq K\}$ and the corresponding part of $\mathbf{s}_{\mathbf{SB}}(t)$ as  $\mathbf{s}_{\mathbf{SB}}(t)[i]$. \\
 Remark that, because $\mathbf{s}$ is $SR$-maintainable  the total number of people in $S_i$ has to remain constant through time. This means that we can view the subsystem given by the seniority-based states $S_{\mathbf{SB}}[i]$ as a Markov chain with at each time step the following three types of flows:
    \begin{itemize}
        \item Incoming flow: the incoming flow concerns transitions from other organisational states into $S_i$ as well as external recruitment into our system via $S_i$. Every member within this category enters $S_i$ in the seniority-based state $S_{0(i)}$.
        \item Staying flow: a part of the population of $S_i$  stays in $S_i$ and obtains additional organisational state seniority. In particular, a staying member moves from $S_{a(i)}$ to $S_{(a+1)(i)}$ in the case $a<K$.
        \item Outgoing flow: the outgoing flow corresponds to transitions from $S_i$ to other organisational states as well as wastage.
    \end{itemize}
    We can use these flows to construct a stochastic matrix. First of all, if we consider the subsystem given by the seniority-based states $S_{\mathbf{SB}}[i]$, it follows that the staying flow is characterised by the corresponding part of $\mathbf{P}_{SM}$, which we denote by $\mathbf{P}_{SM}[i]$. This matrix is, by construction, a strict substochastic matrix. Now, in order to be a $SR$-maintainable structure, every person leaving $S_i$ should be replaced by someone else, who will join the subsystem in the seniority-based state $S_{0(i)}$. Those replacements will either come from other organisational states or from recruitment into $S_{0(i)}$. Note that the fact that the structure $\mathbf{s}$ is $SR$-maintainable implies the existence of a positive recruitment vector in each time step, analogous to the situation in Equation (\ref{maint}).    
    By uniting these types of inflow, we can create a stochastic matrix based on $\mathbf{P}_{SM}[i]$:
    \begin{equation} \label{Pwr2}
        \mathbf{P}[i] := 
        \mathbf{P}_{SM}[i]+
        \begin{pmatrix}
1- \sum_{n=1}^{K+1} \Big( \mathbf{P}_{SM} [i] \Big)_{1 \, n} \\
\vdots \\
1- \sum_{n=1}^{K+1} \Big( \mathbf{P}_{SM} [i] \Big)_{(K+1) \, n} \\
\end{pmatrix} . \begin{pmatrix}
    1, & 0, & 0, & \cdots & 0
\end{pmatrix}
    \end{equation}
    Remark that Equation (\ref{Pwr2}) corresponds to a specific case of Equation (\ref{Pwr}), the case where recruitment is only allowed in $S_{0(i)}$. Furthermore, Equation (\ref{Pwr2})
     allows us to describe the evolution of the population structure related to $S_i$ as a Markov system with states $S_{\mathbf{SB}}[i]$:
    \begin{equation*}
       \mathbf{s}_{\mathbf{SB}}(t)[i]=\mathbf{s}_{\mathbf{SB}}(t-1)[i] \cdot \, \mathbf{P}[i] 
    \end{equation*}
    We claim that the Markov chain with transition matrix $\mathbf{P}[i]$ has a unique fixed point: Because the row sums of $\mathbf{P}_{SM}$ are strictly smaller than $1$, it follows that all elements of the first column of $\mathbf{P}[i]$ are non-zero. Hence, $S_{0(i)}$ can be reached from each state within $S_{\mathbf{SB}}[i]$. We can now partition the set $S_{\mathbf{SB}}[i]$ in two distinct categories of states: either a state belongs to the same communicating class as $S_{0(i)}$ or it does not. If the latter is the case, the number of staff in this state will eventually converge to zero, as there can be no inflow to this state, while the subset of states that do belong to the same communicating class as $S_{0(i)}$ form an ergodic set.  Theorem \, $6.2.1$ in  \citep{kemeny1960finite} states that such a Markov chain has a unique fixed point. \\
    This implies that $((\mathbf{P}[i])^t)_t$ converges:
    \begin{equation*}
        \lim_{t \to \infty} (\mathbf{P}[i])^t ={\mathbf{P}[i]}^*
    \end{equation*}
 and consequently 
    \begin{equation*}
        \lim_{t \to \infty} {\mathbf{s}_{\mathbf{SB}}(t)}[i]=\mathbf{s}_{\mathbf{SB}}(0)[i] \cdot \,{\mathbf{P}[i]}^*:={\mathbf{s}_{\mathbf{SB}}^{*}}[i]
    \end{equation*}
    As this is true for an arbitrary organisational state $S_i$, Equation (\ref{state based structucture converges}) follows.
\end{proof}
In general, there exist different seniority-based paths $(\mathbf{s}_{\mathbf{SB}}(t))_t$ associated to the same $\mathbf{s} \in \mathscr{MR}_{SM}:=\mathscr{MR}_{SM}(1) $. It is easily verified that the property of being associated to the same $\mathbf{s} \in \mathscr{MR}_{SM}$ gives rise to an equivalence relation on the paths. \\
Using Theorem \ref{constant}, it follows that, for each $SR$-maintainable vector $\mathbf{s}$ we can choose the constant path associated with $(\mathbf{s}_{\mathbf{SB}}^*)_t$, satisfying Equations (\ref{SR1}, \ref{SR2}, \ref{SR3}) and where $\mathbf{s}_{\mathbf{SB}}^*$ is maintainable for $\mathbf{P}_{SM}$, as representative of the equivalence class associated to $\mathbf{s}$. Corollary \ref{maintainable region and constant paths} summarizes this insight as a key feature in determining the $SR$-maintainable region $\mathscr{MR}_{SM}$.

\begin{corollary} \label{maintainable region and constant paths}
In order to determine the maintainable region $\mathscr{MR}_{SM}$  for a semi-Markov model with transition matrix $\mathbf{P}_{SM}$, it is sufficient to identify the constant paths $(\mathbf{s}_{\mathbf{SB}}^*)_t$ with $\mathbf{s}_{\mathbf{SB}}^*$ maintainable for $\mathbf{P}_{SM}$. More specifically, the maintainable region $\mathscr{MR}_{SM}$ is the set of population structures $\mathbf{s}_{\mathbf{SB}}^* \mathbf{U}$ with $\mathbf{s}_{\mathbf{SB}}^*$ maintainable for $\mathbf{P}_{SM}$.
\end{corollary}

Remark that Equation (\ref{SR1}), together with the convergence of $(\mathbf{s}_{\mathbf{SB}}(t))_t$ implies the convergence of ${\mathbf{r}}_{\mathbf{SB}}^{+}(t)$. In what follows, the normalisation of the limit recruitment vector is denoted by ${\mathbf{r}}^*$.
The $SR$-maintainable region can be characterised by adapting the procedure to deduce the maintainable region for a Markov chain under control by recruitment as was done in \cite{bartholomew1967stochastic}. As $\mathbf{s}_{\mathbf{SB}}^*$ is a maintainable structure for $\mathbf{P}_{SM}$, we start from Equation (\ref{maineq}) with $\alpha=0$.
\begin{equation}
{\mathbf{s}}_{\mathbf{SB}}^* ={\mathbf{s}}_{\mathbf{SB}}^* \mathbf{P}_{SM} +  {\mathbf{s}}_{\mathbf{SB}}^*{\mathbf{w}}_{\mathbf{SB}}^{'}
{\mathbf{r}}^*
\end{equation}

Note that the substochasticity of $\mathbf{P}_{SM}$ implies the invertibility of $\mathbf{I}-\mathbf{P}_{SM}$, so we can rewrite this equation as follows :
\begin{equation} \label{inversePSM}
{\mathbf{s}}_{\mathbf{SB}}^* =  {\mathbf{s}}_{\mathbf{SB}}^*{\mathbf{w}}_{\mathbf{SB}}^{'}
{\mathbf{r}}^*  (\mathbf{I}-\mathbf{P}_{SM})^{-1}
\end{equation}
If we now multiply both sides of Equation (\ref{inversePSM}) with the re-union matrix $\mathbf{U}$, we obtain:
   \begin{equation} \label{stateM}
{\mathbf{s}}=:{\mathbf{s}}_{\mathbf{SB}}^*  \mathbf{U}={\mathbf{s}}_{\mathbf{SB}}^* {\mathbf{w}_{\mathbf{SB}}^{'}}
{\mathbf{r}^*}   (\mathbf{I}-\mathbf{P}_{SM})^{-1} \mathbf{U}
    \end{equation}
By multiplying the probability vector ${\mathbf{s}}_{\mathbf{SB}}^*$ with $\mathbf{U}$, we obtain the population structure ${\mathbf{s}}$ of length $l$.
If we multiply both sides of Equation (\ref{stateM}) with $E$, the vector of length $l$ with all elements equal to $1$, while decomposing ${\mathbf{r}}^*$ as 
 \begin{equation*}
     {\mathbf{r}^*}
     =r_1 \mathbf{E}_1+ r_2 \mathbf{E}_2 + \cdots  + r_l \mathbf{E}_l
 \end{equation*}
 where $\mathbf{E}_i$ is the standard unit vector of length $(K+1)l$ with the $i$th element equal to $1$, and while writing $d_i$ for the sum of the $i$th row of $(\mathbf{I}-\mathbf{P}_{SM})^{-1}\mathbf{U}$, we obtain the following:
\begin{equation*}
1= {\mathbf{s}}_{\mathbf{SB}}^* {\mathbf{w}_{\mathbf{SB}}^{'}}\Big(\sum_i r_i d_i\Big) \implies {\mathbf{s}}_{\mathbf{SB}}^* {\mathbf{w}_{\mathbf{SB}}^{'}=\frac{1}{\sum_i r_i d_i}}
\end{equation*}
Substituting this result into Equation (\ref{stateM}), we see that
\begin{equation*} 
 {\mathbf{s}}=\frac{1}{\sum_i r_i d_i} {\mathbf{r}^*}   (\mathbf{I}-\mathbf{P}_{SM})^{-1} \mathbf{U} = \sum_j \frac{r_j d_j}{\sum_i r_i d_i} \frac{1}{d_j} \mathbf{E}_j  (\mathbf{I}-\mathbf{P}_{SM})^{-1} \mathbf{U}
\end{equation*}
Writing $\col{\mathbf{R}}^j$ for the $j$th row of $(\mathbf{I}-\mathbf{P}_{SM})^{-1} \mathbf{U}$ we finally obtain:

\begin{equation} \label{SMR2}
  {\mathbf{s}}=\sum_j \frac{r_j d_j}{\sum_i r_i d_i} \frac{1}{d_j} \mathbf{E}_j (\mathbf{I}-\mathbf{P}_{SM})^{-1} \mathbf{U}=\sum_j \frac{r_j d_j}{\sum_i r_i d_i} \frac{\col{\mathbf{R}}^j}{d_j}
\end{equation}
which is a convex combination of the normalised rows of $(\mathbf{I}-\mathbf{P}_{SM})^{-1}\mathbf{U}$.
\begin{remark}
    As recruitment is only allowed in seniority-based states with zero seniority, solely the $r_j$ corresponding to those states are non-zero. I.e. only the rows $\col{\mathbf{R}}^j$ corresponding to the seniority-based states with zero seniority should be taken into account in Equation (\ref{SMR2}).
\end{remark}
So the subsequent Theorem follows:

\begin{theorem} \label{SMM}
    For a semi-Markov model with transition matrix $\mathbf{P}_{SM}$, with constant size and under control by recruitment, the $SR$-maintainable region $\mathscr{MR}_{SM}$ is the convex set with vertices the normalised rows of $(\mathbf{I}-\mathbf{P}_{SM})^{-1}\mathbf{U}$ corresponding to the seniority-based states with zero seniority.
\end{theorem}

\subsection{Extensions to exponential growth and contraction} \label{SRMo}
As is the case for a Markov system (see Theorem \ref{markovexco}), the procedure introduced in Section \ref{SRMcon} can be adapted if we want to retain the relative state sizes in an organisation with constant growth factor $1+\alpha$. \\
First we need the following definition of the matrix $\mathbf{P}[i]^{(1+\alpha)}$ which is inspired by the construction of $\mathbf{P}[i]$ in the proof of Theorem \ref{constant}:
\begin{definition} If we write $\mathbf{P}_{SM}[i]$ for the submatrix of $\mathbf{P}_{SM}$ that corresponds to the flows of the seniority-based states associated with the organisational state $S_i$, we define $\mathbf{P}[i]^{(1+\alpha)}$ in the following way:
    \begin{equation*}
        \mathbf{P}[i]^{(1+\alpha)} := 
        \mathbf{P}_{SM}[i]+
        \begin{pmatrix}
(1+\alpha)- \sum_{n=1}^{K+1} \Big( \mathbf{P}_{SM} [i] \Big)_{1 \, n} \\
\vdots \\
(1+\alpha)- \sum_{n=1}^{K+1} \Big( \mathbf{P}_{SM} [i] \Big)_{(K+1) \, n} \\
\end{pmatrix} . \begin{pmatrix}
    1, & 0, & 0, & \cdots & 0
\end{pmatrix}
    \end{equation*}
\end{definition}
    
We can now state the following generalisation of Theorem \ref{constant}:
\begin{theorem}
    Suppose that we have a matrix $\mathbf{P}_{SM}$ which characterises a semi-Markov model, and the re-union matrix $ \mathbf{U}$. Denote the maximal organisational state seniority by $K$. Suppose that all of the matrices $\mathbf{P}[i]^{(1+\alpha)}$ have a unique fixed point. If the structure $\mathbf{s}$ is $SR$-maintainable for a system with growth factor $1+\alpha$, then, it follows, for a seniority-based path $\mathbf{s}_{\mathbf{SB}}(t)$ associated to $\mathbf{s}$ that the sequence $\mathbf{s}_{\mathbf{SB}}(t)$, converges, i.e
    \begin{equation}\label{state based structucture converges2}
        \lim_{t \to \infty} 
        \mathbf{s}_{\mathbf{SB}}(t) = \mathbf{s}_{\mathbf{SB}}^*
    \end{equation}
\end{theorem}
\begin{proof}
    We can follow the proof of Theorem \ref{constant} with minor adaptations. 
    We will prove the convergence of  $\mathbf{s}_{\mathbf{SB}}(t)$ by demonstrating the convergence of its constituent parts $\mathbf{s}_{\mathbf{SB}}(t)[i]$. As $\mathbf{s}$ is maintainable for a semi-Markov system with growth factor $1+\alpha$ we can again use the incoming, staying and outgoing flows to construct a matrix which encodes the evolution of the subsystem with state space $S_\mathbf{SB}[i]$. 
    Once more, the staying flow is characterised by the corresponding part of $\mathbf{P}_{SM}$, which we denote by  $\mathbf{P}_{SM}[i]$, and this matrix is by definition a strict substochastic matrix. In Theorem \ref{constant}, we had to replace each person leaving the system with another person entering the system, necessarily in the seniority-based state $S_{0(i)}$, a one to one replacement.
    In this case we have to replace each person leaving the system by a \textquotedblleft$(1+\alpha)$\textquotedblright person. More concrete, the row sums of the matrix $\mathbf{P}[i]^{(1+\alpha)}$ should now be equal to $(1+\alpha)$. This leads to
    \begin{equation*}
        \mathbf{P}[i]^{(1+\alpha)} = 
        \mathbf{P}_{SM}[i]+
        \begin{pmatrix}
(1+\alpha)- \sum_{n=1}^{K+1} \Big( \mathbf{P}_{SM} [i] \Big)_{1 \, n} \\
\vdots \\
(1+\alpha)- \sum_{n=1}^{K+1} \Big( \mathbf{P}_{SM} [i] \Big)_{(K+1) \, n} \\
\end{pmatrix} . \begin{pmatrix}
    1, & 0, & 0, & \cdots & 0
\end{pmatrix}
    \end{equation*}
Note that this matrix describes the evolution of the non-normalised population structure related to $S_i$. If we want to consider the $SR$-maintainable vectors, which are probability vectors, we need to consider the stochastic matrix $\frac{\mathbf{P}[i]^{(1+\alpha)}}{1+\alpha}$:
 \begin{equation*}
       \mathbf{s}_{\mathbf{SB}}(t)[i]=\mathbf{s}_{\mathbf{SB}}(t-1)[i] \cdot \, \frac{\mathbf{P}[i]^{(1+\alpha)}}{1+\alpha} 
    \end{equation*}
As the matrix $\frac{\mathbf{P}[i]^{(1+\alpha)}}{1+\alpha}$ has a unique fixed point, it follows that $((\frac{\mathbf{P}[i]^{(1+\alpha)}}{1+\alpha})^t)_t$ converges
 which implies the convergence of 
   $({\mathbf{s}_{\mathbf{SB}}(t)}[i])$.
    Since this is true for arbitrary $S_i$, Equation (\ref{state based structucture converges2}) follows.
\end{proof}
Using the same reasoning that lead to Corollary \ref{maintainable region and constant paths}, we obtain the following:
\begin{corollary} 
In order to determine the maintainable region $\mathscr{MR}_{SM}(1+\alpha)$  for a semi-Markov model with transition matrix $\mathbf{P}_{SM}$ and with growth factor $1+\alpha$, provided that all of the $\mathbf{P}[i]^{(1+\alpha)}$ have a unique fixed point,  
it is sufficient to identify the constant paths $(\mathbf{s}_{\mathbf{SB}}^*)_t$ with $\mathbf{s}_{\mathbf{SB}}^*$ maintainable for $\mathbf{P}_{SM}$ with growth factor $1+\alpha$. 
\end{corollary}
Because $\mathbf{s}_{\mathbf{SB}}^*$ is a maintainable structure for $\mathbf{P}_{SM}$, it follows that $\alpha \geq - \max {w}_i$
and Proposition \ref{a en nonzero column} holds. 
Utilizing the same logic outlined in Section \ref{SRMcon} and coupling it with Theorem \ref{markovexco}, we arrive at the subsequent theorem.
\begin{theorem} \label{smarkovexco}
The maintainable region $\mathscr{MR}_{SM}(1+\alpha)$ of a semi-Markov system with transition matrix $\mathbf{P}_{SM}$ under control by recruitment and with growth factor $1+\alpha$, provided that all of the $\mathbf{P}[i]^{(1+\alpha)}$ have a unique fixed point, satisfies the following:
    \begin{enumerate}
    \item  The maintainable region is a convex set of structures $\mathbf{s}$ that lies in a hyperplane 
    \begin{equation*}
        \mathbf{s}=\mathbf{s}_{\mathbf{SB}}^* \mathbf{U}
    \end{equation*}
 where $\mathbf{s}_{\mathbf{SB}}^*$ is a vector of the intersection of the following halfspaces and hyperplanes:
    \begin{equation} \label{eqhyp}
        \begin{cases}
    {(\mathbf{s}_{\mathbf{SB}}^*  \mathbf{P}_{SM})}_i
    \leq (1+\alpha) ({\mathbf{s}_{\mathbf{SB}}^* })_i & \text{if } \quad i-1   \equiv_{K+1}  0  \\
    {(\mathbf{s}_{\mathbf{SB}}^*  \mathbf{P}_{SM})}_i
    = (1+\alpha) ({\mathbf{s}_{\mathbf{SB}}^* })_i & \text{else } 
\end{cases}
    \end{equation}

    \item If \(  \mathbf{(1+\alpha)I}-\mathbf{P}_{SM} \) is invertible then the maintainable region consists of the  convex combinations of the normalised rows of $(\mathbf{(1+\alpha)I}-\mathbf{P}_{SM})^{-1} \mathbf{U}$ corresponding to the seniority-based  states with zero seniority, that belong to $\Delta^{l-1}$.
    \item If \( \alpha > \rho(\mathbf{P}_{SM}) - 1 \) then the vertices of the maintainable region are the normalised rows of $(\mathbf{(1+\alpha)I}-\mathbf{P}_{SM})^{-1} \mathbf{U}$ corresponding to the seniority-based  states with zero seniority.
\end{enumerate}
\end{theorem}
\begin{proof}
    As the vector $\mathbf{s}_{\mathbf{SB}}^*$ is maintainable for $\mathbf{P}_{SM}$ and non-zero recruitment is only allowed in the seniority based states with zero seniority, Equation (\ref{eqhyp}) follows. Statements $2$ and $3$ follow by the same reasoning that lead to Theorem \ref{SMM}. 
\end{proof}
\newpage
\section{$SR$-maintainable region and maintainable regions of induced Markov chains} \label{section5}
\subsection{Relationship between $SR$-maintainability and maintainability of induced Markov chains}  A semi-Markov chain inherently comprises more information than a Markov chain. Suppose that we have a longitudinal data set at our disposal which admits the estimation of the parameters of a semi-Markov model.  Then we obtain the sequence of matrices $\{\mathbf{P}(k)\}_{0 \leq k \leq K}$ according to Theorem \ref{bigmatrix} and we can use the appropriate re-union matrix $\mathbf{U}$ to study the $SR$-maintainability of this semi-Markov system. Now suppose that we use the same data set to estimate the transition matrix of a Markov model and study the maintainability of that Markov system. The question that arises then is whether there exists any relation between the set of maintainable and the set of $SR$-maintainable population structures.

Remark that it is possible to estimate the transition matrix of a Markov chain starting from $\{\mathbf{P}(k)\}_{0 \leq k \leq K}$ by using for each organisational state $S_i$ the relative frequencies for the seniority-based states in $S_\mathbf{SB}[i]$. \\
 If we supplement the organisational states $\mathbf{S}$ with a wastage state and denote the corresponding transition matrix as $\mathbf{P}_M^S$ we can formulate the  following proposition: 
\begin{proposition}
     Suppose that we have a sequence of $l \times l$ transition matrices \, $\{\mathbf{P}(k)\}_{0 \leq k \leq K}$ \, as in Theorem \ref{bigmatrix} that are estimated on real world data, where the states correspond to the seniority-based states $\mathbf{S_{SB}}$. If we estimate a Markov chain $\mathbf{P}_M^S$ on the same data, we obtain that:
     \begin{equation} \label{columns}
         \mathbf{P}_M^S=
         \sum_{k=0}^{K} \diag(\alpha_k^1,\alpha_k^2, \dots, \alpha_k^l) \cdot \mathbf{P}(k)
    \end{equation}
     Where $\diag(\alpha_k^1,\alpha_k^2, \dots, \alpha_k^l)$ corresponds to an $l \times l$ diagonal matrix with diagonal elements $\alpha_k^1,\alpha_k^2, \dots, \alpha_k^l$, and
    where  $\alpha_k^i$ equals the proportion of the number of people in state $S_i$ with organisational state seniority $k$ relative to the total number of people in the state $S_i$ in the data. 
    
\end{proposition}
   
\begin{proof}
Following the classical notation in \citep{barbu2009semi}, we infer that 
    \begin{equation*}
        \alpha_k^i= \frac{N_i(k)}{N_i}
    \end{equation*}
    The notations $N_{ij}(k)$ and $N_i(k)$ are as in Equation (\ref{schatting P_{ij}(k)}).
    Now we see that:
    \begin{equation*}
        \sum_{k=0}^{K} \alpha_k^i \widehat{P_{ij}(k)}= \sum_{k=0}^{K} \alpha_k^i \frac{N_{ij}(k)}{N_i(k)} = \sum_{k=0}^{K} \frac{N_i(k)}{N_i} \frac{N_{ij}(k)}{N_i(k)}=\frac{1}{N_i} \sum_{k=0}^{K} N_{ij}(k) = \frac{N_{ij}}{N_i}=\left({\mathbf{P}_M^S}\right)_{ij}
    \end{equation*}
    and Equation (\ref{columns}) follows.   
\end{proof}
\begin{remark}
        Note that we used matrix notation to express the fact that the $i$th row of $\mathbf{P}_M^S$ can be seen as a convex combination of the $i$th rows of the matrices $\mathbf{P}(k)$.
    \end{remark} 
So, given a semi-Markov chain as a sequence  of transition matrices \, $\{\mathbf{P}(k)\}_{0 \leq k \leq K}$ and a list of vectors $(\mathbf{\alpha_k})_{{0 \leq k \leq K}}$ we can create the $\alpha-$corresponding Markov model. Observe that, in practice, two data sets resulting in the same (estimated) semi-Markov chain, may result in two different lists of vectors $(\mathbf{\alpha_k})$ and subsequently lead to two different $\alpha-$corresponding Markov chains.  \\
In the next section, we compare the $SR$-maintainable region of the semi-Markov chain $\mathscr{MR}_{SM}$ to the maintainable regions of induced Markov chains $\mathscr{MR}_{M(SM)}$. 

\subsection{Illustration: control theory for a semi-Markov chain with constant size} \label{sec:example}
Suppose that we want to preserve the proportions in the organisational states in a semi-Markov system for which $\mathbf{P}(k)$ is given by: \\
$\mathbf{P}(0)=\left(\begin{matrix}{{0.2}}&{{0.5}}&{{0}}&{{0.3}}\\ {{0}}&{{0.7}}&{{0.2}}&{{0.1}}\\ {{0}}&{{0}}&{{0.9}}&{{0.1}}\\ {{0}}&{{0}}&{{0}}&{{1}}\\ \end{matrix}\right)$, $\mathbf{P}(1)=\left(\begin{matrix}{{0.6}}&{{0.3}}&{{0}}&{{0.1}}\\ {{0}}&{{0.5}}&{{0.45}}&{{0.05}}\\{{0}}&{{0}}&{{0.9}}& {{0.1}}\\ {{0}}&{{0}}&{{0}}&{{1}}\\ \end{matrix}\right)$ and \\
$\mathbf{P}(2)=\left(\begin{matrix}{{0}}&{{0}}&{{0}}&{{1}}\\ {{0}}&{{0}}&{{0}}&{{1}}\\ {{0}}&{{0}}&{{0}}&{{1}}\\ {{0}}&{{0}}&{{0}}&{{1}}\\ \end{matrix}\right)$ \\
we obtain, following Definition \ref{general_P}:  
$$\mathbf{P}_{SM}=$$
\begin{equation} \label{PSM}
\begin{blockarray}{cccccccccc}
&S_{0(1)} & S_{1(1)} & S_{2(1)} & S_{0(2)} & S_{1(2)} & S_{2(2)} & S_{0(3)} & S_{1(3)} & S_{2(3)}  \\
\begin{block}{c(ccccccccc ) @{}}
  S_{0(1)} &\quad \quad 0 & 0.2 & 0 & 0.5 & 0 & 0 & 0 & 0 & 0  \quad\\
  S_{1(1)} &\quad \quad 0 & 0&0.6 & 0.3 & 0 & 0& 0 & 0 & 0 \quad\\
  S_{2(1)} &\quad \quad 0 & 0&0 & 0 & 0 & 0 & 0 & 0 & 0\quad\\
  S_{0(2)} &\quad \quad 0 & 0&0 & 0 & 0.7 & 0 & 0.2& 0& 0 \quad\\
  S_{1(2)} &\quad \quad 0 & 0&0 & 0 & 0 & 0.5 & 0.45& 0& 0  \quad\\
  S_{2(2)} &\quad \quad 0 & 0&0 & 0 & 0 & 0 & 0 & 0 & 0\quad\\
  S_{0(3)} &\quad \quad 0 & 0 & 0 & 0 & 0 & 0 & 0 & 0.9 & 0  \quad\\
S_{1(3)} &\quad \quad 0 & 0 & 0 & 0 & 0 & 0 & 0 & 0 & 0.9  \quad\\
S_{2(3)} &\quad \quad 0 & 0&0 & 0 & 0 & 0 & 0 & 0 & 0\quad\\
\end{block}
\end{blockarray}
\end{equation}
\text{and} ${{\mathbf{w}_{\mathbf{SB}}}}=(0.3,0.1,1,0.1,0.05,1,0.1,0.1,1)$. \\ \\
 Let ${\mathbf{s}}_{\mathbf{SB}}^*=\left(s_{0(1)},s_{1(1)},s_{2(1)},s_{0(2)},s_{1(2)},s_{2(2)},s_{0(3)},s_{1(3)},s_{2(3)}\right)$ 
 and 
 $\mathbf{r}^*=\left(r_{1},0,0,r_{2},0,0,r_{3},0,0\right)$ 
Starting from Equation (\ref{stateM}), we see that:

 \begin{equation*}
{\mathbf{s}}_{\mathbf{SB}}^* U={\mathbf{s}}_{\mathbf{SB}}^* {\mathbf{w}_{\mathbf{SB}}^{'}}
{\mathbf{r}^*}  (\mathbf{I}-P_{SM})^{-1} \mathbf{U}
\end{equation*} with re-union matrix $\mathbf{U}=
\begin{pmatrix}
1 & 0 &0\\
1 & 0 &0\\
1 & 0 &0\\
0 & 1 &0\\
0 & 1 &0\\
0 & 1 &0\\
0&0&1 \\
0&0&1 \\
0&0&1 \\
\end{pmatrix}$. Simplifying this equation yields:\begin{align*}
    {\mathbf{s}}=(s_1,s_2,s_3)&=\left(s_{1,0}+s_{1,1}+s_{1,2},s_{2,0}+s_{2,1}+s_{2,2},s_{3,0}+s_{3,1}+s_{3,2}\right)={\mathbf{s}}_{\mathbf{SB}}^* {\mathbf{w}_{\mathbf{SB}}^{'}}
{\mathbf{r}^*}
(\mathbf{I}-\mathbf{P}_{SM})^{-1} \mathbf{U}
\end{align*}
Which implies that:
\begin{align} \nonumber
  {\mathbf{s}}=&(s_1,s_2,s_3)=  {\mathbf{s}}_{\mathbf{SB}}^* {\mathbf{w}_{\mathbf{SB}}^{'}}
{\mathbf{r}^*} \begin{pmatrix}
1 & \frac{1}{5} & \frac{3}{25} & \frac{14}{25} & \frac{49}{125} & \frac{49}{250} & \frac{721}{2500} & \frac{6489}{25000} & \frac{58401}{250000} \\
0 & 1 & \frac{3}{5} & \frac{3}{10} & \frac{21}{100} & \frac{21}{200} & \frac{309}{2000} & \frac{2781}{20000} & \frac{25029}{200000} \\
0 & 0 & 1 & 0 & 0 & 0 & 0 & 0 & 0 \\
0 & 0 & 0 & 1 & \frac{7}{10} & \frac{7}{20} & \frac{103}{200} & \frac{927}{2000} & \frac{8343}{20000} \\
0 & 0 & 0 & 0 & 1 & \frac{1}{2} & \frac{9}{20} & \frac{81}{200} & \frac{729}{2000} \\
0 & 0 & 0 & 0 & 0 & 1 & 0 & 0 & 0 \\
0 & 0 & 0 & 0 & 0 & 0 & 1 & \frac{9}{10} & \frac{81}{100} \\
0 & 0 & 0 & 0 & 0 & 0 & 0 & 1 & \frac{9}{10} \\
0 & 0 & 0 & 0 & 0 & 0 & 0 & 0 & 1
\end{pmatrix}
\begin{pmatrix}
1 & 0 &0\\
1 & 0 &0\\
1 & 0 &0\\
0 & 1 &0\\
0 & 1 &0\\
0 & 1 &0\\
0&0&1 \\
0&0&1 \\
0&0&1 \\
\end{pmatrix} \\ 
&= {\mathbf{s}}_{\mathbf{SB}}^* {\mathbf{w}_{\mathbf{SB}}^{'}}
\left(r_{1},0,0,r_{2},0,0,r_{3},0,0\right)
\begin{pmatrix}
\frac{33}{25} & \frac{287}{250} & \frac{195391}{250000} \\
\frac{8}{5} & \frac{123}{200} & \frac{83739}{200000} \\
1 & 0 & 0 \\
0 & \frac{41}{20} & \frac{27913}{20000} \\
0 & \frac{3}{2} & \frac{2439}{2000} \\
0 & 1 & 0 \\
0 & 0 & \frac{271}{100} \\
0 & 0 & \frac{19}{10} \\
0 & 0 & 1
\end{pmatrix}\label{exsm}
\end{align}

Multiplying both members of Equation (\ref{exsm}) on the right side with $E=
\begin{pmatrix}
1  \\
1  \\
1  \\
\end{pmatrix}$ yields:
\begin{equation*}
    1={\mathbf{s}}_{\mathbf{SB}}^* {\mathbf{w}_{\mathbf{SB}}^{'}}
 \left(r_1 \frac{812391 }{250000} + r_2 \frac{68913 }{20000} + r_3 \frac{271 }{100}
 \right)
\end{equation*}
So it follows that:
\begin{align*}
     {\mathbf{s}}=(s_1,s_2,s_3)= & \frac{r_1 \frac{812391 }{250000}}{r_1 \frac{812391 }{250000} + r_2 \frac{68913 }{20000} + r_3 \frac{271 }{100}}
     \frac{250000}{812391 }
 \left( \frac{33}{25}, \frac{287}{250}, \frac{195391}{250000} \right)\\
 +&
\frac{r_2 \frac{68913 }{20000}}{r_1 \frac{812391 }{250000} + r_2 \frac{68913 }{20000} + r_3 \frac{271 }{100}}
     \frac{20000}{68913 }
 \left( 0 , \frac{41}{20} , \frac{27913}{20000}  \right)\\
 +&
\frac{ r_3 \frac{271 }{100}}{r_1 \frac{812391 }{250000} + r_2 \frac{68913 }{20000} + r_3 \frac{271 }{100}}
     \frac{100}{271 }
 \left(0 , 0 , \frac{271}{100}  \right) \\
  \approx& a (0.406208, 0.353278, 0.240513)+b(0,
0.594953,
0.405047)+(1-a-b) (0, 0,1)
\end{align*}
With $0 \leq a, b \leq 1$ and $1-a-b \geq 0$. So the maintainable region consists of a convex combination of the vectors $\left(\frac{110000}{270797} ,
\frac{287000}{812391},
\frac{195391}{812391} \right), \left(0 ,
\frac{41000}{68913},
\frac{27913}{68913}\right),$ and $\left(0,0,1\right)$. 
\begin{figure}[h!]
    \centering
    \includegraphics[width=0.8\textwidth]{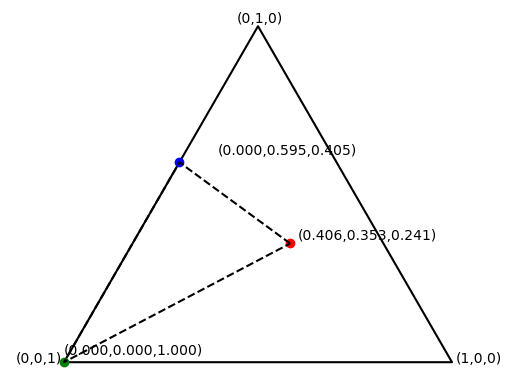}
    \caption{The maintainable region $\mathscr{MR}_{SM}$ for the semi-Markov chain with $\mathbf{P}_{SM}$ as in Equation (\ref{PSM}).}
    \label{fig:maintainable}
\end{figure} \\ \newpage
Reconsider the example above and suppose that $\alpha=(\mathbf{\alpha_0},\mathbf{\alpha_1},\mathbf{\alpha_2})$, where ${\mathbf{\alpha_0}}=\left(0.01,0.28,0.39 \right), {\mathbf{\alpha_1}}=\left(0.68,0.12,0.17\right)$ and ${\mathbf{\alpha_2}}=\left(0.31,0.6,0.44 \right)$. Using this information in combination with the $\mathbf{P}(k)$ as described above, we infer that the $\alpha$-corresponding Markov model is given by:
\begin{equation*}
    \mathbf{P}^I_{M1}=\left(\begin{array}{c c c}{{0.41}}&{{0.209}}&{{0}}\\ {{0}}&{{0.256}}&{{0.11}}\\ {{0}}&{{0}}&{{0.504}}\\ \end{array}\right),
\end{equation*}  \\
Using Theorem \ref{the:mainmarkov}, the maintainable region of this Markov chain, $\mathscr{MR}_{M1(SM)}$, is the convex hull of the normalised rows of  $(\mathbf{I}-\mathbf{P}^I_{M1})^{-1}$. 
\begin{equation*}
    (\mathbf{I}-\mathbf{P}^I_{M1})^{-1}=\begin{pmatrix}
    1.69491525 & 0.47612539 & 0.10559232 \\
    0.0 & 1.34408602 & 0.29808359 \\
    0.0 & 0.0 & 2.01612903 \\
\end{pmatrix}
\end{equation*}
which implies that $\mathscr{MR}_{M1(SM)}$ is the convex hull of the vectors \\
 $(0.7444833137641775, 0.20913576959235633, 0.046380916643466115)$, $(0, \frac{248}{303}, \frac{55}{303})$, and $(0, 0, 1)$.
 Now suppose that $\beta=(\mathbf{\beta_0},\mathbf{\beta_1},\mathbf{\beta_2})$, where ${\mathbf{\beta_0}}=\left(0.15,0.11,0.16 \right), \\{\mathbf{\beta_1}}=\left(0.6,0.56,0.64\right)$ and ${\mathbf{\beta_2}}=\left(0.25,0.33,0.2 \right)$. Using this information we obtain the $\beta$-corresponding Markov model:
\begin{equation*}
    \mathbf{P}^I_{M2}=\left(\begin{array}{c c c}{{0.39}}&{{0.255}}&{{0}}\\ {{0}}&{{0.357}}&{{0.274}}\\ {{0}}&{{0}}&{{0.72}}\\ \end{array}\right),
\end{equation*} \\
The maintainable region for this Markov chain, $\mathscr{MR}_{M2(SM)}$, is given by the convex hull of the normalised rows of  $(\mathbf{I}-\mathbf{P}^I_{M2})^{-1}$. 
\begin{equation*}
    (\mathbf{I}-\mathbf{P}^I_{M2})^{-1}=\begin{pmatrix}
1.63934426 & 0.65012875 & 0.63619742 \\
0.0 & 1.55520995 & 1.52188403 \\
0.0 & 0.0 & 3.57142857
\end{pmatrix}
\end{equation*}
which implies that $\mathscr{MR}_{M2(SM)}$ is the convex hull of the vectors \\
 $(\frac{18004}{32131}, \frac{7140}{32131}, \frac{6987}{32131})$, $(0, \frac{140}{277}, \frac{137}{277})$, and $(0, 0, 1)$. \newline 
\begin{figure}[h!]
    \centering
    \includegraphics[width=0.8\textwidth]{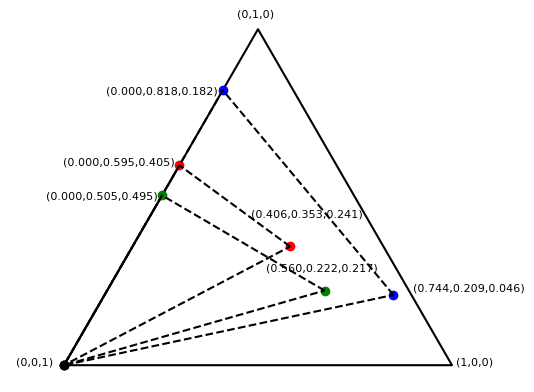}
    \caption{The maintainable regions $\mathscr{MR}_{M1(SM)}$ (blue), $\mathscr{MR}_{M2(SM)}$ (green) and the $SR$-maintainable region (red).}
    \label{fig:maintainable}
\end{figure} \\
We conclude that, in general, there is no obvious connection between the $SR$-maintainable region $\mathscr{MR}_{SM}$ of a semi-Markov chain and the maintainable region $\mathscr{MR}_{M(SM)}$ of an induced Markov chain. While for some induced Markov chains, as for $M=M1$, $\mathscr{MR}_{SM} \subseteq \mathscr{MR}_{M(SM)}$ holds, this is not the case for others such as for $M=M2$. Even the opposite inclusion might hold, which is for example the case for the $\gamma$-corresponding model, where $\gamma=(\mathbf{\gamma_0},\mathbf{\gamma_1},\mathbf{\gamma_2})$, with${\mathbf{\gamma_0}}=\left(0.31,0.29,0.04 \right), {\mathbf{\gamma_1}}=\left(0.27,0.48,0.57\right)$ and ${\mathbf{\gamma_2}}=\left(0.42,0.23,0.39 \right)$.

\section{Conclusions and further research avenues}
In this paper, we investigate control by recruitment, extending the traditional framework of maintainability in Markov chains to accommodate for semi-Markov chains. The focus of our investigation is not merely an extension but also a novel approach by introducing a new type of maintainability termed as State Re-union Maintainability ($SR$-maintainability). Furthermore, we obtain essential theoretical insights that enable the determination of the $SR$-maintainable region.
Upon analysis, we find that direct comparisons between the maintainable regions of Markov chains and the $SR$-maintainable region of an associated semi-Markov chain yield no generalisable conclusions regarding the relationship between them.  \\
The strength of our research lies in its innovative approach to generalising maintainability and introducing $SR$-maintainability, offering both theoretical insights and an algorithmic solution for practitioners, who can determine the $SR$-maintainable region by the use of Theorem \ref{smarkovexco}.
\\
Future lines of inquiry might delve into the potential generalisations of Re-union matrices. This exploration would involve extending the concept of $SR$-maintainability to encompass the maintainability of various combinations of seniority-based states, such as regrouping according to total seniority or salary grade. Such an extension could enhance our comprehension of how diverse internal states can be merged, potentially opening doors to the exploration of partial ($SR$-) maintainability. This, in turn, would permit the utilization of different $\mathbf{U}$-matrices, potentially broadening the scope of applicability of the $SR$-maintainability concept.
Hence, our study serves as a stepping stone in the landscape of maintainability in stochastic processes, combining traditional Markov chains with their more flexible semi-Markov counterparts, and introducing a new theoretical construct in the form of $SR$-maintainability.
\pagebreak
\newpage
\bibliography{references}
\end{document}